\documentclass[10pt,reqno]{amsart}
\usepackage{amssymb,enumerate,color,fancyvrb,prooftree,bm}
\usepackage[colorlinks=true,urlcolor=blue,pagebackref=true]{hyperref}
\usepackage{multicol}
\usepackage{setspace}
\allowdisplaybreaks
\theoremstyle{plain}
\newtheorem{theorem}{Theorem}[section]

\newtheorem{lemma}[theorem]{Lemma}
\newtheorem{corollary}[theorem]{Corollary}

\theoremstyle{definition}
\newtheorem{definition}[theorem]{Definition}
\newtheorem{notation}[theorem]{Notation}
\theoremstyle{remark}
\newtheorem{remark}[theorem]{Remark}
\newtheorem{problem}[theorem]{Problem}

\newcommand\tb{\textbf}

\newcommand{\F}{\mathbb{F}}
\newcommand{\Gbb}{\mathbb{G}}
\newcommand{\Ff}{\mathbf{F}}
\newcommand{\Gf}{\mathbf{G}}
\newcommand\ft{\Phi}
\newcommand\ar{P}
\newcommand\Gr{\Gbb(\ar)}
\newcommand\mb{\mathbf}      
\newcommand\mt{\mathtt}
\newcommand{\vp}{\varphi}
\newcommand{\ha}{{\sf HA}} 
\newcommand{\czf}{{\sf CZF}}
\newcommand{\izf}{{\sf IZF}}
\newcommand{\acft}{{\sf AC}_{\mathrm{FT}}}
\newcommand{\ac}{{\sf AC}}
\newcommand{\markov}{{\sf MP}}
\newcommand{\ip}{{\sf IP}}
\newcommand{\ipr}{{\sf IPR}}
\newcommand{\dc}{{\sf DC}}
\newcommand{\rdc}{{\sf RDC}}
\newcommand{\pax}{{\sf PAx}}
\newcommand{\rea}{{\sf REA}}
\newcommand{\fo}{\Vdash_{tr}}   
\newcommand{\ap}{\mb p}         
\newcommand{\pair}[1]{\langle{#1}\rangle}
\newcommand{\va}{\mathrm{V}_{tr}(A)}
\newcommand{\V}{\mathrm{V}}
\newcommand{\vset}[1]{\{{#1}\}_{A}}
\newcommand{\vpair}[1]{\langle{#1}\rangle_{\!A}}

\newcommand{\vo}[1]{{#1}^\circ}
\newcommand{\vs}[1]{{#1}^*}
\DeclareMathOperator{\dom}{\mathrm{dom}}
\DeclareMathOperator{\op}{\mathrm{op}}
\DeclareMathOperator{\fun}{\mathrm{Fun}}

\DeclareMathOperator{\ps}{\mathcal{P}}
\DeclareMathOperator{\succe}{\mb{succ}}
\DeclareMathOperator{\pred}{\mb{pred}}
\newcommand{\imp}{\rightarrow}

\newcommand{\biimp}{\leftrightarrow}
\newcommand{\Biimp}{\Leftrightarrow}

\def\rul#1#2#3{\prooftree     #1 \justifies #2 \using{#3} \endprooftree}

\title[Choice and independence of premise rules in set theory]{Choice and independence of premise rules in intuitionistic set theory}
\author[Frittaion]{Emanuele Frittaion}
\address{Department of Mathematics, Technische Universit\"{a}t Darmstadt, Germany}
\email{frittaion@mathematik.tu-darmstadt.de}
\thanks{Emanuele Frittaion's research was supported by the Alexander von Humboldt Foundation.}
\author[Nemoto]{Takako Nemoto}
\address{Graduate School of Information Sciences, Tohoku University,
	6-3-09 Aoba, Aramaki-aza Aoba-ku, Sendai, 980-8579, Japan}
\email{nemototakako@gmail.com}
\thanks{Takako Nemoto's research was supported by the Japan Society for the Promotion of Science
(JSPS), Core-to-Core Program (A. Advanced Research Networks).}
\author[Rathjen]{Michael Rathjen}
\address{Department  of Pure Mathematics,
	University of Leeds, UK}
\email{M.Rathjen@leeds.ac.uk}
\thanks{Michael Rathjen's research was supported by the John Templeton
Foundation (A new dawn of intuitionism: mathematical and philosophical advances, ID
60842).}

\begin{document}
	
\subjclass[2020]{Primary: 03F03; Secondary: 03F25, 03F50, 03F55}
	
\begin{abstract}
Choice and	independence of premise principles play an important role in characterizing Kreisel's modified realizability and  G\"odel's Dialectica interpretation. In this paper we show that a great many intuitionistic set theories are closed under the corresponding rules for finite types   over $\mathbb N$. It is also shown that the existence property (or existential definability property) holds for statements of the form 
$\exists y^{\sigma}\, \vp(y)$, where the variable $y$ ranges over objects of finite type $\sigma$. This applies in particular to $\czf$ (Constructive Zermelo-Fraenkel set theory)  and $\izf$ (Intuitionistic Zermelo-Fraenkel set theory), two systems known not to have the general existence property.
On the technical side, the paper uses a method that amalgamates generic realizability for set theory with truth,  whereby the underlying partial combinatory algebra is required to {\em contain}  all objects of finite type. 
\end{abstract}
	
\maketitle

\section{Introduction}

There are (at least) three types of  classically valid principles  that figure prominently in constructive mathematics:  $\acft$ ({\em Choice in Finite Types}), $\markov$ ({\em Markov's principle}) and  $\ip$ ({\em Independence of Premise principle}).  All three are required for a characterization of G\"odel's Dialectica interpretation (see \cite[3.5.10]{T73}, \cite[Proposition 8.13]{K08}),
whereas Kreisel's modified realizability for intuitionistic finite type arithmetic $\ha^{\omega}$ is axiomatized by $\acft$ and $\ip$ alone. To be more precise, let
\begin{eqnarray*}
\ac_{\mathrm{FT}} && \forall x^{\sigma}\, \exists y^{\tau}\, \vp(x,y)\to \exists f^{\sigma\tau}\, \forall x^{\sigma}\, \vp(x,fx) \\[1ex]
\ip_{\neg} && (\neg\psi\imp \exists y^\sigma\, \vp(y))\imp \exists y^\sigma\, (\neg\psi\imp \vp(y)) \\[1ex]
\ip_{\mathrm{ef}} && (\psi_{\mathrm{ef}}\to \exists y^{\sigma}\, \vp(y))\to \exists y^{\sigma}\, (\psi_{\mathrm{ef}} \to \vp(y))
\end{eqnarray*}
where $\sigma,\tau$ signify  finite types, $z^{\rho}$ varies over objects of finite type $\rho$, and $\psi_{\mathrm{ef}}$ is assumed to be $\exists$-free, i.e., it neither contains existential quantifiers nor disjunctions.\footnote{Of course, it is also assumed that $y$ is not a free variable of $\neg\psi$ and $\psi_{\mathrm{ef}}$.}

Then the following holds (see e.g.\ 	\cite[3.4.8]{T73}, \cite[Theorem 3.7]{T98}, \cite[Theorem 5.12]{K08}):
\begin{theorem}  With  $\Vdash_{mr}$ signifying modified realizability, we have:
\begin{itemize} 
\item[(i)] $\ha^\omega+\acft+\ip_{\mathrm{ef}}\vdash \vp\biimp \exists x\,(x\Vdash_{mr} \vp)$.
\end{itemize}
Let $T\in \{\ha^\omega, {\sf E\text{-}HA}^\omega\}$.\footnote{For a definition of $\ha^\omega$ and its extensional variant ${\sf E\text{-}HA}^\omega$ cf.\ e.g.\ \cite{T73} or \cite{K08}.} Then:
\begin{itemize}
\item[(ii)]   $T +\acft+\ip_{\mathrm{ef}}\vdash \vp$ iff $T \vdash (t\Vdash_{mr} \vp)$ for some term $t$.
\end{itemize}
\end{theorem}
An important application  of modified realizability (to be correct, its truth variant, aka modified realizability with truth) is that (i) $T$ is closed under the  rule of choice in finite types, (ii)  $T$ is closed under the independence of premise rule
for negated  formulas $\mathsf{IPR}_{\neg}$ (in particular, it is closed under $\ipr_{\mathrm{ef}}$), and (iii)  $T$ satisfies explicit definability. Short and sweet (see e.g.\ \cite[Theorem 3.8]{T98}, \cite[Corollary 5.24]{K08}):
\begin{theorem}\label{haupt}
Let $T\in\{{\sf HA}^\omega, {\sf E\text{-}HA}^\omega\}$. Then:
\begin{itemize}
\item[(i)] If $T\vdash \forall x^{\sigma}\, \exists y^{\tau}\, \vp(x,y)$, then $T\vdash \exists f^{\sigma\tau}\, \forall x^{\sigma}\, \vp(x,fx)$.
\item[(ii)] If $T\vdash \neg\psi\to \exists y^{\sigma}\, \vp(y)$, then $T\vdash \exists y^{\sigma}\, (\neg\psi \to \vp(y))$;
\item[(iii)] If $T\vdash \exists y^{\sigma}\, \vp(y)$, then $T\vdash \vp(t)$, for a suitable term $t$.
\end{itemize}
\end{theorem}

The current paper shows that results similar to Theorem \ref{haupt} hold for a great  many set theories $T$, including 
$\czf$ (Constructive Zermelo-Fraenkel set theory) and $\izf$ (Intuitionistic Zermelo-Fraenkel set theory),  even if augmented by suitable choice principles and large set axioms.\footnote{The intuitionistic rendering of large cardinal axioms.} A more precise delineation of the kind of set theories eligible for this theorem is that $T$ should be self-validating with respect to generic 
realizability combined with truth (cf.\ Theorem \ref{sound}).

\begin{theorem}[see Theorems \ref{choice rule}, \ref{premise rule} and \ref{def}]\label{Main} 
An array of set theories $T$ including $\czf$ and $\izf$ satisfy the following:
\begin{itemize}
\item[(i)] If $T\vdash \forall x^{\sigma}\, \exists y^{\tau}\, \vp(x,y)$, then $T\vdash \exists f^{\sigma\tau}\, \forall x^{\sigma}\, \vp(x,fx)$;
\item[(ii)] If $T\vdash \forall x\, (\neg \psi(x) \to \exists y^\sigma\, \vp(x,y))$, then 
		$T\vdash \exists y\, \forall x\, (\neg \psi(x) \to  y\in \sigma \land \vp(x,y))$;
\item[(iii)] If $T\vdash \forall x\, (\forall z \, \vartheta(x,z) \to \exists y^\sigma\, \vp(x,y))$ and $T\vdash \forall z\, (\vartheta(x,z)\,\vee\,\neg \vartheta(x,z))$, then 
		$T\vdash \exists y^\sigma\, \forall x\,(\forall z\, \vartheta(x,z) \to \vp(x,y))$;
\item[(iv)] If $T\vdash \exists y^\sigma\, \vp(y)$, then $T\vdash \exists ! y^\sigma\, (\delta(y)\land \vp(y))$,   for some formula $\delta(y)$.
\end{itemize}
\end{theorem} 

To properly shelve the perhaps perplexing results, i.e.,  from a classical viewpoint,   it is good to bear in mind that our theorem applies to set theories
$T$  closed under the following rules:
\begin{itemize}
	\item[(v)] the {\em Unzerlegbarkeits rule}, namely, if $T\vdash \forall x\, (\vp(x)\lor \psi(x))$, then $T\vdash \forall x\, \vp(x)$ or $T\vdash\forall x\, \psi(x)$;
	\item[(vi)] the {\em Uniformity rule}, namely,  if $T\vdash \forall x\, \exists y\in\omega\, \vp(x,y)$, then $T\vdash \exists y\in\omega\, \forall x\, \vp(x,y)$
\end{itemize}
(see \cite[Theorem 1.2]{R05}  and \cite[Theorem 7.4]{R08}).\footnote{Discussions of uniformity under realizability seem to appear first in the literature in 
Friedman's \cite{F73} and Troelstra's \cite{T73a}.} A by-product of Theorem \ref{Main} (cf.\ Section \ref{uniformity}) is that  rule (vi)  still holds  when  $y$ ranges over objects of a given finite type.  

It is known that  $\izf$ (see \cite{My73}) and $\czf$ (see \cite{R05}) have the numerical existence property, so Theorem \ref{Main} part (iv) extends this property to a larger collection of existential formulas. On the other hand, it is known by
work of Friedman and \v{S}\v{c}edrov \cite{FS85} that $\izf$ does not have  the  general existence property $\mathsf{EP}$   while Swan \cite{S14} proved that $\mathsf{EP}$ also fails for $\czf$. In the latter case the culprit is $\czf$'s  Subset Collection axiom  as shown by the third author \cite{Rathjen2023} since
the version of $\czf$ with Exponentiation in lieu of Subset Collection  has the $\mathsf{EP}$.

On the technical side, the paper uses the method of generic realizability with truth from \cite{R05}.  The main novelty of the paper, however, is the introduction   of a special kind of  partial combinatory algebra (pca) to the realizability framework. We single out  the notion of a {\em reflexive pca over finite types}   to deal specifically with objects of finite type.   An instance of this new phenomenon, and its construction in weak set theories such as $\czf$, is duly supplied.

Generic realizability is better understood as a generalization of Kleene's 1945 realizability rather than Kreisel's modified realizability. Even so, it differs from both.  On one hand,  generic realizability (with truth) combined with a reflexive pca over finite types bears close similarities to modified realizability (with truth). Roughly, a truth realizer of a $\forall x^\sigma\, \exists y^\tau$ statement yields a choice functional of type $\sigma\tau$. This gives closure under the rule of choice in all finite types. On the other hand, totality is a crucial aspect of modified realizability with truth in establishing closure under independence of premise rule for (extensional) finite type arithmetic:  realizers are total functionals.  As it turns out,  the independence of premise rules of Theorem \ref{Main} do not require the use of a total combinatory algebra:   the inherent nature of set theory (a set is given by its members) allows us to dispense with totality.  Note in contrast that truth variants of Kleene's realizability (e.g.\ $q$-realizability \cite{T73, TD88}) do not yield closure under independence of premise rule for first order arithmetic $\ha$, and one has to resort to other methods such as Kleene's or Aczel's slash  instead.

Our approach to establish independence of premise rules however is not without shortcomings, as can be seen from Theorem \ref{Main} part (ii), and we do not know whether the following more genuine version  holds true.

\begin{problem}
	Is the following an admissible rule of $\czf$ or any other familiar constructive/intuitionistic set theory $T$?
	\begin{itemize}
		\item If $T\vdash \neg\psi\to \exists y^\sigma \, \vp(y)$,
		then $T\vdash\exists y^\sigma \, (\neg\psi\to\vp(y))$,
		where  $y$ is not free in $\psi$.
	\end{itemize}
\end{problem}
Notice that  as a special case of  Theorem \ref{Main} part (ii) one only obtains the following weaker rule: 
\begin{itemize}
	\item if $T\vdash \neg\psi\to \exists y^\sigma \, \vp(y)$,
	then $T\vdash\exists y\, (\neg\psi\to y\in\sigma\land \vp(y))$,
\end{itemize} 
where the $y$ in the conclusion is guaranteed to be of type $\sigma$ only in case the  premise $\neg\psi$ is verified.

More in general, closure under independence of premise rule with no type restrictions  remains an open problem: 
\begin{problem}
Is the following  an admissible rule of $\czf$ or any other familiar constructive/intuitionistic set theory $T$?
\begin{itemize}
	\item	If $T\vdash \neg\psi\to \exists y\, \vp(y)$,
		then $T\vdash\exists y\, (\neg\psi\to\vp(y))$,
		where  $y$ is not free in $\psi$.
\end{itemize}
\end{problem}
The present  paper shows that this holds true when $\exists y$ is bounded by some finite type. 

\subsection{Generic realizability for set theory}
Realizability semantics are ubiquitous in the study of intuitionistic theories. In the case of set theory, they differ in important aspects from Kleene's \cite{K45} realizability in their treatment of the quantifiers. Its origin is Kreisel's and Troelstra's \cite{KT70} definition of realizability for second order Heyting arithmetic.
This was applied to systems of higher order arithmetic and (intensional) set theory by
Friedman \cite{F73} and Beeson \cite{B85}. McCarty \cite{M84} and \cite{M86} adapted Kreisel-Troelstra realizability directly to  extensional intuitionistic set theories such as $\izf$.  This type of
realizability can also be formalized in $\czf$ (see \cite{R06}) to yield a
self-validating semantics for $\czf$. 

Realizability combined with truth appears in connection with function realizability in Kleene \cite{K69} and was also studied by others (see \cite{T98} for the history). Troelstra considers realizability with truth in the arithmetic context and  in connection with 
modified realizability in \cite[1.6, 2.1, 3.4]{T98}. 
In generic realizability for extensional set theory, however, the background universe $\V$ and the realizability universe $\V(A)$ erected over a partial combinatory algebra $A$ are rather different ``worlds", and it is prima facie not clear how to view a statement as  talking about a state of affairs in $\V$ and  $\V(A)$ at the same time.
The  paper \cite{R05} introduced a new realizability structure
$\va$ that arises by amalgamating the realizability structure  $\V(A)$
with the universe of sets in a coherent, albeit  rather
complicated way. This approach to realizability with truth based on  $\va$ will also be used in the present paper.\footnote{For more information, the introduction of \cite{R05} contains 
a historical account of realizability for set theories and the roots of generic realizability in particular.} A very rough heuristics for using this method of realizability with truth is that it often works with principles that are validated in a  realizability model  based on a particular partial combinatory algebra and that switching to the corresponding version with the additional truth component one can derive the pertaining rule.\footnote{E.g. for Church's rule and Troelstra's Uniformity rule this was done in \cite{R05} using Kleene's first algebra.} 

\subsection{Comparison with other approaches to showing the independence of premise rule}
Closure under $\mathsf{IPR}$  has been shown  for $\sf{HAS}$ (second order  Heyting arithmetic)  by Troelstra \cite[2.10]{T73a} and for $\sf{HAH}$ (Heyting  arithmetic in higher types, aka intuitionistic type theory) by Lambek and Scott \cite{LS81}.\footnote{For a definition of $\sf{HAS}$ and $\sf{HAH}$ cf.\ \cite[pp.\ 164 and 170]{TD88}. $\sf HAS$ is a subsystem of $\czf$ + (Full Separation). The two theories are known to be equiconsistent as shown by Lubarsky \cite{Lubarsky2006}.  $\sf HAH$ is a  fragment of intuitionistic Kripke-Platek set theory plus Powerset, $\sf IKP$ + (Powerset), but strictly weaker in terms of proof-theoretic strength. The latter theory is much weaker than intuitionistic Power Kripke-Platek set theory (see \cite{Rathjen2014}), $\mathsf{IKP}(\mathcal{P})$, which proof-theoretically equates to a version of the Calculus of Constructions with one universe by \cite[Theorem 15.1]{Rathjen2012a}.
And $\mathsf{IKP}(\mathcal{P})$ is just a small fragment of $\mathsf{IZF}$.}
Its admissibility is often established as a by-product of the existence property $\mathsf{EP}$ (also called the existential definability property) that such systems enjoy. 
As to the methods used in the metamathematics of $\sf{HAS}$ and $\sf{HAH}$, one can roughly group them as follows:
\begin{enumerate} 
	\item Proof-theoretic methods: study of the proof structure either in natural deduction systems or sequent calculi (e.g. Prawitz \cite{prawitz70,prawitz71}, Scarpellini \cite{scarpellini70,scarpellini71}, Troelstra \cite{T73a}, Hayashi \cite{hayashi77,hayashi78,hayashi80}).
	\item Functional interpretation (e.g. Girard \cite{girard71,girard73}).
	\item Extensions of Kleene's slash method \cite{K69a} (e.g.  Moschovakis \cite{M67}, Myhill \cite{My73,My75}, Friedman \cite{F73}, Lambek-Scott \cite{LS80}).
	\item Topos-theoretic methods, conceptualizing term models as topoi (``free topos'') and using techniques such as  Freyd covers and topos glueing (e.g. Freyd \cite{freyd},  Moerdijk \cite{moerdijk}, Lambek-Scott
	\cite{LS81,LS83,LS86},
	\v{S}\v{c}edrov-Scott \cite{SS82}).
\end{enumerate}
Taking intuitionistic Zermelo-Fraenkel set theory, $\izf$, to be the ``typical'' set theory of this paper, one can perhaps immediately say that such strong theories are currently not amenable to the methods of  (1) and (2).  
The topos-theoretic methods of (4) have turned out to be equivalent in a strong sense to Friedman's modification  of Kleene's slash in \cite{F73}: ``{\em Thus Freyd's use of retracts and Friedman's impredicative assignment of indices turn out to be one and the same process}'' \cite[443]{SS82}.
In view of the foregoing, we will only consider the Friedman slash method and point out the obstacles one faces when attempting to obtain  the results of this paper via this method. 
For any of the slash methods to apply to a theory $T$ one needs a language that has sufficiently many terms to serve as names for the objects that are describable in  $T$. In the case of arithmetic this is easy as one has the numerals. In the context of higher order systems or set theories one is usually compelled to move to a theory $T^*$ with a richer language that
is conservative over $T$. Typically let $T$ be a set theory with explicit set existence axioms, i.e. with axioms that define the contents of the set being asserted to exists, namely if it is of the form \begin{eqnarray}\label{explicit} &&\forall x_0\ldots \forall x_{n-1}\,[\psi(x_0,\ldots , x_{n-1}) \to \exists y\,\forall u\,[u\in y \leftrightarrow \vp(u,x_0,\ldots,  x_{n-1})].\end{eqnarray} One then simultaneously inductively defines a new theory $T^*$  and a collection of closed terms ${\mathcal T}^*$  such that $T^*$ comprises the axioms of $T$ 
and for each axiom of the form (\ref{explicit}), whenever there is $t_0,\ldots,t_{n-1}\in {\mathcal T}^*$  with $T^*\vdash \psi(t_0,\ldots,t_{n-1})$, one adds a new constant $c_{\vp,\psi}(t_0,\ldots,t_{n-1})$ to $ {\mathcal T}^*$ and an axiom 
\[\forall u\,[u\in c_{\vp,\psi}(t_0,\ldots,t_{n-1})\leftrightarrow \vp(u,t_0,\ldots,t_{n-1})]. \]
It turns out that $T^*$ is conservative over $T$. In a further step one defines the system $T^{\diamond}$ that is obtained by splitting up each term of $T^*$ into many terms. Thus the terms of $T^{\diamond}$ are of the form $c_{\vp,\psi}(t_0,\ldots,t_{n-1})^X$, where $X$ is any set of closed  terms of $T^*$ satisfying certain conditions. Echoing Myhill's words \cite[p.\ 369]{My75},  roughly $c_{\vp,\psi}(t_0,\ldots,t_{n-1})^X$ denotes the set 
$\{u\mid \vp(u,t_0,\ldots,  t_{n-1})\}$ and $X$ is the ``reason'' that we know  $\vp(u,t_0,\ldots, t_{n-1})$. For our purposes it is not  necessary to spell out the details.
It suffices to know that the clauses for the quantifiers in the definition of the Friedman slash refer to the terms of $T^{\diamond}$. This slash interpretation works for many theories with explicit set 
existence axioms such as intuitionistic Zermelo-Fraenkel set theory when based on Replacement rather than Collection. As a typical application one obtains the set existence property.
However, the Friedman slash is not known to work for theories whose set existence axioms are not explicit such as  the Collection, Strong Collection, Subset Collection, the Regular Extension Axiom and the Presentation Axiom. Indeed,  as shown by Friedman and \v{S}\v{c}edrov \cite{FS85},  $\izf$ does not have the $\mathsf{EP}$ and Swan \cite{S14} proved that $\czf$ also lacks  the $\mathsf{EP}$, rendering it  unlikely that the Friedman slash can be applied to these theories to establish closure under  independence of premise rules whereas the method of realizability with truth works perfectly well.

\section{Intuitionistic set theory} 

The language of constructive Zermelo-Fraenkel set theory $\czf$ is same first order
language as that of classical Zermelo-Fraenkel set theory $\sf ZF$ whose only
non-logical symbol is the binary predicate $\in$.
We use $x$, $y$, $z$, $u$, $v$, $w$, possibly with superscripts, for variables in the
language of $\czf$. 
The logic of $\czf$ is intuitionistic first order logic with equality.
The axioms of $\czf$ are as follows (universal closures):

\noindent\tb{Extensionality:} $\forall z\, (z\in x\leftrightarrow z\in y)\to x=y$.

\noindent\tb{Pairing:} $\exists z\, (x\in z\land y\in z)$.

\noindent\tb{Union:} $\exists y\, \forall z\, ( \exists w\in x\, (z\in w)\to z\in y)$.

\noindent\tb{Infinity:} $\exists w\, \forall x\, [x\in w\leftrightarrow \forall y\, (y\notin
	x)\lor \exists y\in w\, \forall z\, (z\in x\leftrightarrow
	z\in y\lor z=y)]$.
 
\noindent\tb{Set Induction:} 
	$\forall x\, [ \forall y\in x\, \varphi(y)\to \varphi(x)]\to \forall
	x\, \varphi(x)$, for any formula $\varphi$.
	
\noindent\tb{Bounded Separation:} $ \exists y\, \forall z\, [z\in y\leftrightarrow z\in x\land
	\varphi(z)]$, for any {\it bounded} formula $\varphi$. A formula is {\it bounded} or
	{\it restricted} if it is constructed from prime formulae using $\to$, $\neg$,
	$\land$, $\lor$, $\forall x\in y$ and $\exists x\in y$ only. 
	
\noindent\tb{Strong Collection:} 
	\begin{align*}
		\forall y\in x\, \exists z\, \varphi(y,z)\to 
		\exists w\, [\forall y\in x\, \exists z\in w\, \varphi(y,z)\land \forall z\in
		w\, \exists y\in x\, \varphi(y,z)],
	\end{align*}
for any formula $\varphi$.

\noindent\tb{Subset Collection:}
	\begin{multline*}
		\forall x\, \forall y\, \exists z\, \forall u\, [\forall v\in x\, \exists w\in
		y\, \varphi(v,w,u)\to \\
		\exists y'\in z\, [\forall v\in x\, \exists w\in y'\, \varphi(v,w,u)\land \forall
		w\in y'\, \exists v\in x\, \varphi(v,w,u)]],
	\end{multline*}
for any formula $\varphi$.

In what follows, we shall assume that the language $\czf$ has constants $\emptyset$
denoting the {\it empty set}, 
$\omega$ denoting the set of von Neumann natural numbers.
One can take the axioms $\forall x\, (x\notin\emptyset)$ for $\emptyset$
and $\forall x\, [x\in\omega\leftrightarrow (x=\emptyset\lor \exists y\in\omega\, \forall
z\, (z\in x\leftrightarrow z\in y\lor z=y))]$ for $\omega$.
We write $x+1$ for $x\cup\{ x\}$ and use $n$, $m$, and $l$ for elements of $\omega$.

We consider also several extensions of $\czf$ with other principles.

\noindent\tb{Full Separation:} $\exists y\, \forall z\, [z\in y\leftrightarrow z\in x\land
	\varphi(z)]$, for any formula $\varphi$.

\noindent\tb{Powerset:} $ \exists y\, \forall z\, (z\subseteq x\to z\in y)$.

The system $\czf+(\mathbf{Full\;Separation})+(\mathbf{Powerset})$ is called $\izf$
(cf. \cite{My73} or \cite[VIII.1]{B85}). 

\noindent\tb{Markov principle} ($\markov$)\tb{:}   If $\forall
	n\in\omega\, (\varphi(n)\lor\neg\varphi(n))$ and 
	$\neg\neg\exists n\in\omega\, \varphi(n)$, then $\exists n\in\omega\, \varphi(n)$.

\noindent\tb{Axiom of Countable Choice} ($\ac_{\omega}$)\tb{:} If  $\forall n\in\omega\, \exists x\, \vp(n,x)$, then $\exists f\, (f \text{ is a function }\land \dom(f)=\omega\land \forall n\in\omega\, \vp(n,f(n)))$. 
	
\noindent\tb{Dependent Choices Axiom} ($\dc$)\tb{:} If $\forall x\in z\, \exists y\in z\, \vp(x,y)$, then $\forall x\in z\, \exists f\, (f \text{ is a function }\land \dom(f)=\omega\land f(0)=x\land \forall n\in\omega\, \vp(f(n),f(n+1)))$.

\noindent\tb{Relativised Dependent Choices Axiom} ($\rdc$)\tb{:} If $\forall x\, (\psi(x)\imp \exists y\, (\psi(y)\land \vp(x,y))$, then for every $x$ such that $\psi(x)$ there is a function $f$ such that  $\dom(f)=\omega\land f(0)=x\land \forall n\in\omega\, \vp(f(n),f(n+1))$.

\begin{definition}
A set $x$ is {\em inhabited} if $\exists y\, (y\in x)$.
An inhabited set $x$ is {\it regular} if $x$ is transitive, and for every $y\in x$ and a
	set $z\subseteq y\times x$ if $\forall u\in y\, \exists v\, (\langle u,v\rangle\in z)$, then
	there is a set $w\in x$ such that
	\begin{align*}
		\forall u\in y\, \exists v\in w\, (\langle u,v\rangle\in z)\land \forall v\in w\, \forall u\in
		y\, (\langle u,v\rangle\in z).
	\end{align*}
\end{definition}

\noindent\tb{Regular Extension Axiom}  ($\rea$)\tb{:} Every set is a subset of a regular set.

\begin{definition}
	A set $x$ is {\it projective} if for any $x$-indexed family $(y_u)_{u\in x}$ of inhabited
	sets $y_u$, there exists a function $f$ with domain $x$ such that $f(u)\in y_u$ for all $u\in x$.
\end{definition}

\noindent   {\bf Presentation Axiom}  ($\pax$)\tb{:} Every
set is the surjective image of a projective set.

\subsection{Finite types}

\begin{definition}
Finite types $\sigma$ and their extensions $F_\sigma$ are defined by the following clauses:
\begin{itemize}
		\item $o\in \ft$ and $F_o=\omega$;
		\item if $\sigma,\tau\in\ft$, then $(\sigma)\tau\in\ft$ and $F_{(\sigma)\tau}={F_\sigma} \to {F_\tau}=\{ \text{total functions from $F_\sigma$ to $F_\tau$}\}$.
\end{itemize} 
If there is no risk of confusion, we  write $\sigma\tau$ or $\sigma\to\tau$ to denote the type $(\sigma)\tau$. 
The sets $\ft$ of all finite types,  $\mb F=\{ F_\sigma\colon \sigma\in\ft\}$ and $\F=\bigcup\Ff$
all exist in $\czf$.
\end{definition}

\begin{definition}\label{theta}
	There is a formula $\vartheta(\sigma,z)$ (also written $\vartheta_{\sigma}(z)$) such that:
	\begin{itemize}
		\item $\vartheta(o,z)\biimp z \text{ is } \omega$;
		\item $\vartheta(\sigma\tau,z)\biimp \forall z_\sigma\, \forall z_\tau\, (\vartheta(\sigma,z_\sigma)\land \vartheta(\tau,z_\tau)\imp \forall f\, (f\in z\biimp  \fun(f,z_\sigma,z_\tau)))$.
	\end{itemize}
	Here, $\fun(f,x,y)$ is an abbreviation for \lq\lq $f$ is a function from $x$ to $y$."
\end{definition}

\begin{notation}[Official]
	$\forall x^\sigma\ldots$ stands for $\forall z_\sigma\, (\vartheta_\sigma(z_\sigma)\imp \forall x\in z_\sigma\ldots)$. Similarly for $\exists x^\sigma$.
\end{notation}

\section{Partial combinatory algebras I}

In order to define a realizability interpretation we must have a notion of realizing functions to
hand. A particularly general and elegant approach to realizability builds on structures which
have been variably called partial combinatory algebras, applicative structures, or Sch\"onfinkel
algebras. For more information on these structures see \cite{F75,F79, B85,Oosten08}.

\begin{definition} $(A,\cdot)$ is said to be a {\em partial algebra} if $A$ is a set and $\cdot$ is a binary function on some subset of $A\times A$.
Since $\cdot$ is only partial on $A$,  it is convenient to talk about application terms of $(A,\cdot)$, where these terms might not denote an object in $A$.  Given an infinite collection  $x,y,z,\ldots$ of variables, the inductive definition of application terms   is as follows: every variable $x$ and every $a\in A$ is an application 
term; if $s,t$ are application terms then $(s\cdot t)$ is an application term. A closed application term is one without variables and it denotes an element $a$ of $A$ iff it is  $a$  itself or else it is of the form   $(s\cdot t)$ and there are
$b,c\in A$ such that $s$ denotes $b$, $t$ denotes $c$, $(b,c)$ is in the domain of $\cdot$ and  $a=b\cdot c$. If a closed application term $t$ denotes, then we convey this by writing $t\downarrow$.

For application terms $s_1,\ldots,s_r$ we shall just write $s_1\ldots s_r$ to refer to  the application term inductively defined by letting $s_1\ldots s_{n+1}$ be $((s_1\ldots s_n)\cdot s_{n+1})$; 
 so the convention is to drop $\cdot$ and
assume the bracketing to be arranged to the left. We also use $s=t$ to convey that the closed application terms $s$ and $t$ denote the same object in $A$; in particular $s=t$ entails that $s\downarrow$ and $t\downarrow$. We also introduce the  very helpful abbreviation $t\simeq s$ for $(t\downarrow\,\vee\,s\downarrow) \to s=t$.

A {\em partial  combinatory algebra}  (pca) is a partial algebra $(A,\cdot)$  such that $A$ has at least two elements and there are elements 
	$\mb k$ and $\mb s$ in $A$ such that $\mb k a$, $\mb sa$ and  $\mb sab$ are always defined  and 
	\begin{itemize}
		\item $\mb ka b\simeq a$;
		\item $\mb sabc\simeq ac(bc)$
	\end{itemize}
	holds for all $a,b,c\in A$.
\end{definition}
The combinators $\mb k$ and $\mb s$ are due   to Sch\"onfinkel \cite{S24}   while the axiomatic treatment, although formulated just in the total case,  is due to Curry \cite{C30}.  Employing the axioms for the combinators $\mb k$ and $\mb
s$ one can deduce an abstraction lemma yielding $\lambda$-terms (cf.\ \cite{F75}).
\begin{lemma}[abstraction lemma]
For each application term $t(x_1,\ldots,x_{n+1})$, there
is a closed application term $a$, denoted $\lambda x_1\cdots x_n. t$, such that  for all $a_1,\ldots,a_n,b\in A$
\begin{itemize}
\item $aa_1\cdots a_n \!\downarrow$;
\item $aa_1\cdots a_nb\simeq t(a_1,\ldots,a_nb)$. 
\end{itemize}
\end{lemma}

The most important consequence of the Abstraction Lemma is the
Recursion Theorem. It can be derived in the same way as for the
$\lambda$--calculus (cf. \cite{F75}, \cite{F79}, \cite[VI.2.7]{B85}). 

\begin{corollary}[recursion theorem] For every $n>0$ there exists a closed application term $f$ such that for  all $a,b_1,\ldots,b_n\in A$ we have:
\begin{itemize}
	\item $fa\downarrow$; 
    \item  $fab_1\cdots b_n\simeq a(fa) b_1\cdots b_n$. 
	\end{itemize}
\end{corollary}

In every pca, one has pairing and unpairing\footnote{Let $\mb p=\lambda xyz.zxy$, $\mb{p_0}:=\lambda x.x\mb k$, and $\mb{p_1}:=\lambda x.x\bar{\mb k}$, where $\bar{\mb k}:=\lambda xy.y$. Projections $\mb{p_0}$ and $\mb{p_1}$ need not be total. } combinators $\mb p$,  $\mb {p_0}$, and  $\mb {p_1}$ such that: 
\begin{itemize}
	\item  $\mb pab\downarrow$;
	\item   $\mb {p_i}(\mb pa_0a_1)\simeq a_i$.
\end{itemize}

The notion of a pca  is slightly impoverished  compared to that of a model of  Beeson's theory ${\sf PCA}^+$ 
\cite[VI.2]{B85} or Feferman's theory for applicative structures $\sf APP$ (\cite{F75,F79}, \cite[9.3]{TD88}). Although, as Curry showed, every pca can be expanded to a model of $\mathbf{PCA}^+$, which at the same time is also an applicative structure (see \cite[VI.2.9]{B85}), we spell out the sort of structure we are interested in, that is, a model of ${\sf PCA}^+ $. Details 
will become more pertinent 
when we engineer specific ones that  include all finite types (Definition \ref{pca F}).   

\begin{definition}\label{pca omega}
	We say that $A$ is a {\em pca over} $\omega$ if there are extra combinators $\succe$, $\pred$ (successor and predecessor combinators), $\mb d$ (definition by cases combinator), and a  map
	$n\mapsto \bar n$ from $\omega$ to $A$ such that 
	for all $n,m\in \omega$ and $a,b\in A$
	\begin{align*}
		\succe \bar n&\simeq \overline{n+1}, & \pred\overline{n+1}&\simeq \bar n, &
\mb d\bar n\bar mab\simeq
	\begin{cases} a & n=m;\\ b & n\neq m. \end{cases} \end{align*}
	One then defines $\mb 0:=\bar 0$ and $\mb 1:=\bar 1$. 
\end{definition}

\section{Realizability with truth}

\subsection{The general realizability structure}

\cite{R05} introduces a  realizability structure with truth over Kleene's first algebra.
In this paper, we define it over an arbitrary set-sized pca $A$ (both $A$ and the graph $\{ (x, y, z)\in A^3: xy\simeq z\}$ are sets).  

\begin{notation}
	For an ordered pair $x=\pair{x_0,x_1}$, let 
	\begin{align*}
		\vo x&= x_0\\
		\vs x&= x_1.
	\end{align*}	
\end{notation} 

\begin{definition}\label{grs}
Ordinals are transitive sets whose elements are transitive also. As per usual, we use lower case Greek letters $\alpha$ and $\beta$ to range over ordinals. Let $A$ be a pca.
	Besides $\V_{\alpha}$ and $\V$,
	we define $\va_\alpha$ and $\va$ as follows:
\begin{align}
\label{eg} \va_\alpha &
		= \bigcup_{\beta\in\alpha}\bigl \{\pair{x,\hat x} \mid x\in\V_{\beta}\land \hat x\subseteq A\times\va_\beta\land 
		\forall \pair{a,u}\in \hat x\, (u^\circ \in x)\bigr\}
		\\ \nonumber \V_{\alpha}&= \bigcup_{\beta\in
			\alpha}\ps(\V_{\beta}) \\
		\nonumber \va &=
		\bigcup_{\alpha} \va_\alpha\\
		\nonumber \V&= \bigcup_{\alpha} \V_{\alpha}. 
	\end{align}
\end{definition}
As the power set operation  is not available in $\czf$, it is not clear whether the classes $\V$ and $\va$ can be formalized in $\czf$. However, employing the fact that $\czf$ accommodates
inductively defined classes, the classes $\V_\alpha$ and $\va_\alpha$ can be defined in the same vein as in \cite[Lemma 3.4]{R06}.

\begin{lemma}[$\czf$]\label{stets} The following holds:
	\begin{itemize}
		\item[(i)] {\em $\V$ and $\va$} are cumulative: for
		$\beta\in\alpha$, {\em $\V_{\beta}\subseteq\V_{\alpha}$ and
			$\va_\beta\subseteq \va_\alpha$.}
		\item[(ii)] For all sets $x$,  {\em $x\in \V$}.
		\item[(iii)]  If $x,\hat x$ are sets,
		{\em $\hat x\subseteq A\times\va$} and
		$\forall \pair{a,u}\in \hat x\, (u^\circ \in x)$,
		then {\em $\pair{x,\hat x}\in\va$}.
	\end{itemize}
\end{lemma}
\begin{proof}
	This is proved in the same way as \cite[Lemma 4.2]{R05}.
\end{proof}

The definition of $\va_\alpha$ in (\ref{eg}) is perhaps a bit involved. Note first that all the elements of $\va$ are ordered pairs $\pair{x,\hat x}$ such that $\hat x\subseteq A\times\va$. For an ordered pair $\pair{x,\hat x}$ to enter $\va_\alpha$ the first conditions to be met are that $x\in\V_{\beta}$ and $\hat x\subseteq A\times\va_\beta$ for some $\beta\in\alpha$.
Furthermore, it is required that  enough elements of $x$ live in the transitive closure of $\hat x$ in that whenever $\pair{a,u}\in \hat x$ then $\vo u \in x$.

For all intents and purposes,  the following equivalent definition of $\va$ is perfectly justifiable in $\czf$. 

\begin{definition}[universe]
	Given a pca $A$, we inductively define the class $\va$ by the following clause:  
	\begin{itemize}
		\item if $\hat x\subseteq A\times \va$ and for every $\pair{a,u}\in \hat x$ we have $u^\circ\in x$, then $\pair{x,\hat x}\in\va$.
	\end{itemize}
\end{definition}

\begin{definition}[canonical name]
	Let
	\[ \check{x}=\pair{x,\{\pair{\mb 0, \check{u}}\colon u\in x\}}.\]
	Then $\check{x}\in\va$ and $\check{x}^\circ=x$.
\end{definition}

\subsection{Defining realizability}

We now proceed to define a notion of realizability with truth over $\va$, where $A$ is any pca over $\omega$. 

\begin{definition}
Given a formula $\vp$ with parameters in $\va$, let $\vp^\circ$ be the formula obtained by replacing each parameter $x$ in $\vp$ with $\vo x$.
\end{definition}

\begin{notation} We use $(a)_i$ or simply $a_i$ for $\mb {p_i}a$.   Whenever we write an application term $t$, we assume that it is defined. In other words, a formula $\vp(t)$ stands for $\exists a\, (t\simeq a\land \vp(a))$.
\end{notation}

The following truth variant of generic realizability is  due to Rathjen (see \cite{R05}). Bounded quantifiers are treated as quantifiers in their own right, i.e., bounded
and unbounded quantifiers are treated as syntactically different kinds of quantifiers. The
subscript in $\fo$ is supposed to serve as a reminder of \lq\lq realizability with truth."

\begin{definition}[realizability with truth] We define the relation $a\fo \vp$, where $a\in A$ and $\vp$ is a formula with parameters in $\va$. The atomic cases are defined by transfinite recursion. 
	\begin{align*}   
		a&\fo x\in y && \Biimp&& \vo x\in \vo y \land  \exists z\, (\langle a_0,z\rangle \in \vs y\land a_1\fo x=z)\\
		a& \fo x=y &&  \Biimp && \vo x=\vo y  \land \forall \langle b,z\rangle \in \vs x\, ((ab)_0\fo z\in y)  \\
		&&&&& \quad\quad\quad\ \land  \forall \langle b,z\rangle \in \vs y\, ((ab)_1\fo z\in x)\\
		a& \fo \vp\land \psi && \Biimp&& a_0\fo \vp \land a_1\fo \psi \\
		a& \fo \vp\lor\psi &&  \Biimp&& a_0\simeq \mb 0\land a_1\fo \vp \text{ or }  
		  a_0\simeq \mb 1\land a_1\fo \psi \\
		a&\fo \neg\vp && \Biimp&& \neg \vp^\circ\land \forall b\, \neg (b\fo \vp) \\
		a&\fo \vp\imp\psi && \Biimp&& \vp^\circ\imp \psi^\circ\land \forall b\fo \vp\,   (ab\fo \psi) \\
		a& \fo \forall x\in y\, \vp && \Biimp&& \forall x\in \vo y\, \vo\vp\land  \forall \langle b,x\rangle\in \vs y\, (ab\fo \vp) \\
		a&\fo \exists x\in y\, \vp && \Biimp&& \exists x\, (\langle a_0,x\rangle \in \vs y\land a_1\fo \vp)\\
		a& \fo \forall x\, \vp && \Biimp&& \forall x\in \va\, (a\fo \vp) \\
		a& \fo \exists x\, \vp && \Biimp&& \exists x\in \va\, (a\fo \vp)
	\end{align*}
\end{definition}

\begin{lemma} \label{rt}
$\czf$ proves
\[(a\fo \vp)\imp \vp^\circ. \]	
\end{lemma}

\begin{proof}
By induction on the build up of $\vp$. The case of an unbounded universal quantifier follows from the fact that every set has a name in $\va$.
\end{proof}

\begin{lemma}
Negated formulas are self-realizing, that is to say, $\czf$ proves 
\[ \neg \vp^\circ \imp (\mb 0\fo \neg \vp) \biimp \forall a\, (a \fo \neg\vp). \]
\end{lemma}
\begin{proof}
Assume $\neg\vp^\circ$. From $a\fo \vp$, we would get $\vp^\circ$  by Lemma \ref{rt}. But this is absurd.
Hence $\forall a\, \neg (a\fo \vp)$,  and therefore $\mb 0\fo \neg\vp$. The second part is similar.
\end{proof}

\begin{theorem}[Soundness] \label{sound}
Let $T$ be any combination of $\czf$ with
the axioms and schemes (Full Separation), (Powerset), $\sf REA$, 
$\markov$, $\ac_{\omega}$, $\dc$, $\rdc$, and
$\sf PAx$.
Then, for every theorem $\theta$ of $T$, there exists an application term $\mb t$ such that
$T\vdash (t\fo \theta)$. In particular, $\czf$, $\czf+ {\sf REA}$, $\izf$, $\izf+{\sf REA}$
satisfy this property.
Moreover, the proof is effective in that the application term
$\mb t$ can be constructed from the $T$-proof of $\theta$.
\end{theorem}
\begin{proof}
This is proved in the same way as \cite[Theorem 6.1, Theorem 7.2]{R05} and
\cite[Theorem 7.4]{R08}.
\end{proof}

\begin{notation} 
We write $\fo \vp$ for $\exists a\in A\, (a\fo \vp)$.
\end{notation}

\subsection{Pairing}

\begin{definition}[internal pairing]
For $x,y\in\va$, let
\[  \vset{x}=\pair{\{x^\circ\},\{\pair{\mb 0,x}\}}, \]
\[   \vset{x,y}=\pair{\{x^\circ,y^\circ\},\{\pair{\mb 0,x},\pair{\mb 1,y}\}}, \]
\[ \vpair{x,y}=\pair{\pair{x^\circ,y^\circ},\{\pair{\mb 0,\vset{x}}, \pair{\mb 1,\vset{x,y}}\}}. \]

Note that all these sets are in $\va$.
\end{definition}

\begin{notation}
To avoid confusion, let $\op(z,x,y)$ be a formula expressing that $z$ is the ordered pair of $x$ and $y$, that is, $z=\pair{x,y}=\{\{x\},\{x,y\}\}$.	
\end{notation}

As expected, all the desired properties of pairing are realized. Below we list some.
\begin{lemma}\label{pairs}
There are closed application terms $\mb v$, $\mb w$, $\mb z$ such that for all $x,y,z,u,v\in \va$
\begin{align*}
\mb v &\fo \op(\vpair{x,y},x,y), \\
\mb w& \fo \vpair{x,y}=\vpair{u,v}\imp x=u\land y=v,\\
\mb z&\fo \op(z,x,y) \imp z=\vpair{x,y}.
\end{align*}
\end{lemma}

\section{Partial combinatory algebras II}

To deal with the rules of choice and independence of premise in all finite types, we will use our truth variant of generic realizability with pca's {\em containing} all objects of finite type.

\begin{definition}[pca over $\F$]\label{pca F}
We say that $A$ is a {\em pca over} $\F$ if there are  extra combinators $\succe$, $\pred$ (successor and predecessor combinators), $\mb d$ (definition by cases combinator), and a one-to-one map $x\mapsto \bar x$ from $\F$ to $A$ such that
\begin{itemize}
	\item $\bar f \bar  x\simeq \overline{f(x)}$ for $f\in F_{\sigma\tau}$ and $x\in F_\sigma$;
	\item for all $n,m\in \omega$ and $a,b\in A$
	\begin{align*}
		\succe \bar n&\simeq \overline{n+1}, & \pred\overline{n+1}&\simeq \bar n, &
		\mb d\bar n\bar mab\simeq
		\begin{cases} a & n=m;\\ b & n\neq m. \end{cases} \end{align*}	
\end{itemize}

\end{definition}

The idea is to have \emph{nice} names of the form
\[  \dot F_\sigma=\pair{F_\sigma,\{\pair{\bar x,\dot x}\mid x\in F_\sigma\}}\]
for every type $\sigma$. Indeed, we require a little bit more.

\begin{definition}[reflexive pca over $\F$]\label{fte}
A pca $A$ over $\F$ is  \emph{reflexive} on $\F$ if for all $\sigma$ and $\tau$ there is a combinator $\mb i_{\sigma\tau}$ such that 
	\[  \mb i_{\sigma\tau} a\simeq \bar f, \]
     whenever $f\in F_{\sigma\tau}$ and $a\bar x=\overline{f(x)}$ for every $x\in F_\sigma$. 
\end{definition}

Unless otherwise stated, from now on we posit a reflexive pca $A$ over $\F$ within $\czf$. In Section \ref{example} we will give an example and show how to carry out such construction in $\czf$.

\section{Injective names for finite type objects}

\begin{definition}[canonical names for objects of finite type and extensions]
Let $A$ be a pca over $\F$. Let
\[ \dot{\omega}=\pair{\omega,\{\pair{\bar n,\dot n}\mid n\in\omega\}}, \]
where
\[ \dot n=\pair{n,\{\pair{\bar m, \dot m}\mid m<n\}}. \]

For higher types, let
\[ \dot F_{\sigma\tau}=\pair{F_{\sigma\tau},\{ \pair{\bar f,\dot f}\mid f\in F_{\sigma\tau}\}}, \]
where 
\[   \dot f=\pair{f,\{\pair{\bar x, \vpair{\dot x,\dot y}}\mid x\in F_\sigma\land f(x)=y\}}. \]
\end{definition}

\begin{lemma}
Let $A$ be a pca over $\F$. For every $\sigma$,
\begin{itemize}
	\item  $(\dot{x})^\circ= x$ for every $x\in F_\sigma$,
	\item  $(\dot F_\sigma)^\circ= F_\sigma$.
\end{itemize}
\end{lemma}

\begin{definition}[injective]\label{def inj}
A name $x\in\va$ is \emph{injective}  if 
\begin{enumerate}[(i)]
	\item $x^\circ=\{u^\circ\mid \exists a\in A\, (\pair{a,u}\in x^*)\}$;
	\item if $\pair{a,u},\pair{b,v}\in x^*$, then  $a=b$ iff $u^\circ=v^\circ$.
\end{enumerate}
In other words, $\{\pair{a,u^\circ}\mid \pair{a,u}\in x^*\}$ is one-to-one function 
from  $\{a\in A\mid \exists u\, (\pair{a,u}\in x^*)\}$ onto 
$x^\circ$. We say that $x^\circ$ is \emph{injectively presented}. 
\end{definition}

\begin{lemma}\label{inj}
Let $x\in\va$ be injective. Then
\begin{align*}
	  a\fo \forall u\in x\, \vp(u)   && \Biimp&& \forall \pair{b,u}\in x^*\, ab\fo \vp(u).
\end{align*}	
\end{lemma}
\begin{proof}
Indeed, condition (i) is sufficient.
\end{proof}

\begin{lemma}\label{dot}
Let $A$ be a pca over $\F$. For every $\sigma$,
\begin{itemize}
	\item $\dot x$ is injective for every $x\in F_\sigma$;
	\item $\dot F_\sigma$ is injective;
	\item $\fo \dot x=\dot y$ implies $x=y$ for all $x,y\in F_\sigma$ (absoluteness).
\end{itemize} 
\end{lemma}

\section{Realizing finite types}

We want to show that $\czf$ proves $\fo \vartheta_\sigma(\dot F_\sigma)$ for every $\sigma$, provided that $A$ is a reflexive pca over $\F$. Recall that $\vartheta_\sigma(z)$ is the formula asserting that $z$ is the set of all objects of type $\sigma$ (Definition \ref{theta}).

\begin{theorem}[natural numbers]
There exists  a closed application term $\mb e$ such that $\czf$ proves
\[ \mb e \fo \vartheta_o(\dot \omega). \]
\end{theorem}
\begin{proof}
See \cite[Theorem 6.1 (Infinity)]{R05}.	Note that any pca would do the job.
\end{proof}

For arrow types, we use reflexivity.

\begin{theorem}[arrow types]\label{arrow}
For all finite types $\sigma$ and $\tau$ there exists a closed application term $\mb e$ such that $\czf$ proves
	\[   \mb e\fo \forall f\, (f\in \dot F_{\sigma\tau}\biimp \fun(f,\dot F_\sigma,\dot F_\tau)). \]
\end{theorem}
\begin{proof}
Fix types $\sigma$ and $\tau$. It suffices to look for closed application terms $\mb a$ and $\mb b$ such that
\[  \mb a\fo \forall f\in \dot F_{\sigma\tau}\, \fun(f,\dot F_\sigma,\dot F_\tau), \]
	and for every $g\in \va$,
\[  \mb b\fo \fun(g,\dot F_\sigma,\dot F_\tau) \imp g\in \dot F_{\sigma\tau}. \]
	
To ease notation, we identify $x$ with $\bar x$. Since $(\dot F_\sigma)^\circ=F_\sigma$ for every type $\sigma$,  we just need to verify the second half of the pertaining clauses, namely, 
\begin{itemize}
	\item for every $f\in F_{\sigma\tau}$,  $\mb af\fo \fun(\dot f,\dot F_\sigma,\dot F_\tau)$;
	\item if $a\fo \fun(g,\dot F_\sigma, \dot F_\tau)$, then $\mb ba\fo  g\in \dot F_{\sigma\tau}$.
\end{itemize}

For $\mb a$,  we need to find $\mb {r}, \mb {t}, \mb {f}$ such that for every $f\in F_{\sigma\tau}$,
	\[  \tag{1} \mb {r} f\fo \forall z\in \dot f\, \exists x\in\dot F_\sigma\, \exists y\in\dot F_\tau\, \op(z,x,y), \]
	\[   \tag{2 }\mb {t} f\fo \forall x\in\dot F_\sigma\, \exists y\in\dot F_\tau\, \exists z\in \dot f\, \op(z,x,y), \]
	\[  \tag{3} \mb{f} f\fo \forall z_0\in \dot f\, \forall z_1\in\dot f\, \forall x\, \forall y_0\, \forall y_1\, (\op(z_0,x,y_0)\land \op(z_1,x,y_1)\imp y_0=y_1). \]
	
For (1), let \[ \mb {r}=\lambda ax.\mb p x(\mb p (ax)\mb v), \]
where $\mb v\fo \op(\vpair{x,y},x,y)$ for all $x,y\in \va$. Let us verify that $\mb {r}$ does the job. Let $f\in F_{\sigma\tau}$. We want to show that
	\[    \mb {r} f=\lambda x.\ap  x(\mb p (fx)\mb v)\fo \forall z\in\dot f\, \exists x\in\dot F_\sigma\, \exists y\in\dot F_\tau\, \op(z,x,y). \]
As $\dot f$ is injective, by Lemma \ref{inj} it suffices to show that for $x\in F_\sigma$ and $y=f(x)\in F_\tau$ we have that
	\[   \ap x(\mb p (fx)\mb v)\fo \exists x_0\in \dot F_\sigma\, \exists y_0\in\dot F_\tau\, \op(\vpair{\dot x,\dot y},x_0,y_0). \]
By definition, $\pair{x,\dot x}\in(\dot F_\sigma)^*$. Let us check that
	\[ \mb p (fx)\mb v\fo  \exists y_0\in\dot F_\tau\, \op(\vpair{\dot x,\dot y},\dot x,y_0). \]
Note that $fx\simeq y$ and $\pair{y,\dot y}\in(\dot F_\tau)^*$. 
Finally, 
	\[ \mb v\fo \op(\vpair{\dot x,\dot y},\dot x,\dot y). \]
	
For (2), use
	\[   \mb {t}=\lambda ax.\ap (ax)(\ap x\mb v), \]
where $\mb v$ is as above. 
	
For (3), use the fact that $\dot{x}^\circ = x$ for every $x\in F_\sigma$ and the properties of pairing and equality.   

We now construct  $\mb b$.  Here is where the  $\sigma\tau$ combinator $\mb i_{\sigma\tau}$ comes into play. Suppose that
\[   a\fo  \fun(g,\dot F_\sigma, \dot F_\tau). \]
We aim for 
\[ \mb b a\fo g\in \dot F_{\sigma\tau}. \]
We have 
\[  \tag{4} a_0\fo \forall z\in g\, \exists x\in\dot F_\sigma\, \exists y\in\dot F_\tau\, \op(z,x,y), \]
\[ \tag{5} a_{10}\fo \forall x\in\dot F_\sigma\, \exists y\in\dot F_\tau\, \exists z\in g\, \op(z,x,y), \]
\[  \tag{6} a_{11}\fo \forall z_0\in g\, \forall z_1\in g\, \forall x\, \forall y_0\, \forall y_1\, (\op(z_0,x,y_0)\land \op(z_1,x,y_1)\imp y_0=y_1). \]

Since $g^\circ\in (\dot F_{\sigma\tau})^\circ=F_{\sigma\tau}$,  we only need to find $\mb b$ such that 
\[  (\mb b a)_0\simeq f\in F_{\sigma\tau} \text{ and } (\mb b a)_1\fo g=\dot f. \]
	
It follows from (5) that for every $x\in F_\sigma$ there exists (a unique) $y\in F_\tau$ such that $(a_{10}x)_0\simeq y$ and
\[ (a_{10}x)_1\fo  \exists z\in g\, \op(z,\dot x,\dot y). \]
We now apply the  $\sigma\tau$ combinator $\mb i_{\sigma\tau}$. We then have $\mb i_{\sigma\tau}\lambda x.(a_{10}x)_0 \simeq f$, for the (unique) $f\in F_{\sigma\tau}$ such that 
	\[   (a_{10}x)_0\simeq f(x), \text{ for all } x\in F_\sigma. \]
	
We set
	\[ \mb {b}=\lambda a.\ap (\mb i_{\sigma\tau} \lambda x.(a_{10}x)_0)(\mb ha), \]	

and we are left to find $\mb h$ such that
	\[ \mb ha\fo g=\dot f. \]
	
By using (5), it is not difficult to show that $g^\circ(x)=f(x)$ for every $x\in F_\sigma$, and hence $g^\circ=f=(\dot f)^\circ$. It remains to prove the second half of the pertaining clause.

($\subseteq$) Let $\pair{b,z}\in g^*$. We aim for $(\mb hab)_0\fo z\in\dot f$. By (4), there are $x\in F_\sigma$ and $y_0\in F_\tau$ such that
	\[   (a_0b)_0\simeq x\text{ and } (a_0b)_{10}\simeq y_0 \text{ and } (a_0b)_{11}\fo \op(z,\dot x,\dot y_0).\]
On the other hand, by (5), there is a $\pair{b_1,z_1}\in g^*$ with $(a_{10}x)_{10}\simeq b_1$ such that
	\[  (a_{10} x)_{11}\fo \op(z_1,\dot x,\dot{y}), \]
where $y=f(x)$. By (6),
\[ a_{11}bb_1(\ap (a_0b)_{11} (a_{10} x)_{11})\fo \dot y_0=\dot y. \] 


By absoluteness (Lemma \ref{dot})  it follows from $\fo \dot y_0=\dot y$  that $y_0=y$. Then $\mb h$ such that 
\[  (\mb h ab)_0\simeq \ap (a_0b)_0(\mb q(a_0b)_{11}), \]
where $\mb q$ is some fixed term such that 
\[ \mb q\fo \op(z,x,y) \imp z=\vpair{x,y}, \]
is as desired.
	
($\supseteq$) Let $\pair{ x,\vpair{\dot x,\dot y}}\in (\dot f)^*$. By definition, $f(x)=y$.	We aim for $(\mb hax)_1\fo \vpair{\dot x,\dot y}\in g$. Just let
\[ (\mb hax)_1=\ap (a_{10}x)_{10}(\mb r (a_{10}x)_{11}), \]
	where $\mb r$ is some fixed term such that 
	\[   \mb r\fo \op(z,x,y)\imp \vpair{x,y}=z. \]
\end{proof}

\begin{theorem}\label{type}
For every finite type $\sigma$  there exists a closed application term $\mb f$ such that $\czf$ proves
\[ \mb f \fo \vartheta_\sigma(\dot F_\sigma). \]
\end{theorem}

\section{Admissible rules}

\subsection{Choice}

\begin{lemma}[choice for injective  names]\label{choice inj}
	$\czf$ proves
	\[ (\fo \forall u\in x\, \exists v\in y\, \vp(x,y))\imp \exists f\colon x^\circ\to y^\circ\, \forall u\in x^\circ\, \vp^\circ(u,f(u)), \] 
	for all injective names $x,y\in\va$.
\end{lemma}
\begin{proof}
	Suppose 
	\[    e\fo \forall u\in x\, \exists v\in y\, \vp(u,v). \]
	
	Unraveling the definition, we have  that if $\pair{a,u}\in x^*$ then $ea\fo \exists v\in y\, \vp(u, y)$, and hence there is $\pair{b,v}\in y^*$, where $(ea)_0\simeq b$, such that $(ea)_1\fo \vp(u,v)$. In particular,  $\vp^\circ(u^\circ,v^\circ)$. Let
	\[  f=\{\pair{u^\circ,v^\circ}\colon \pair{a,u}\in x^* \land \pair{b,v}\in y^*\land (ea)_0\simeq b\}. \]
	Then $f$ is as desired. In fact, $\dom(f)=x^\circ$ follows from Definition \ref{def inj}\, (i). The fact that $f$ is indeed a function follows from Definition \ref{def inj}\, (ii) applied to both $x$ and $y$. 
\end{proof}


\begin{theorem}[choice rule]\label{choice rule}
	$\czf$ is closed under 
	\[ \rul{\forall x^\sigma\, \exists y^\tau\, \vp(x,y)}{\exists f^{\sigma\tau}\, \forall x^\sigma\, \vp(x,f(x))}{} \]
	Same for $\izf$ and any other theory from Theorem \ref{sound}.
\end{theorem}
\begin{proof}
Use a reflexive pca over $\F$. By soundness, let $\mb e\fo \forall x^\sigma\, \exists y^\tau\, \vp(x,y)$. By Theorem \ref{type} and soundness, we can compute $a$ such that
\[ a \fo \forall x\in \dot F_\sigma\, \exists y\in \dot F_\tau\, \vp(x,y). \]
By injectivity and  Lemma \ref{choice inj}, we conclude 
\[ \exists f\colon F_\sigma\to F_\tau\, \forall x\in F_\sigma\, \vp(x,f(x)), \]
that is, $\exists f^{\sigma\tau}\, \forall x^\sigma\, \vp(x,f(x))$.
\end{proof}

\subsection{Uniformity rules}\label{uniformity}

\begin{theorem}
$\czf$ is closed under
\[ \tag{{\sf UZR}} \rul{\forall \bm x\, (\vp(\bm x)\lor \psi(\bm x))}{\forall \bm x\, \vp(\bm x)\quad \text{or} \quad  \forall \bm x\, \psi(\bm x)}{} \]
Same for $\izf$ and any other theory from Theorem \ref{sound}.
\end{theorem}
\begin{proof}
Use generic realizability with truth and Kleene's first algebra. See \cite[Theorem 1.2]{R05}.
\end{proof}

\begin{theorem}
$\czf$ is closed under 
\[ \tag{${\sf UR}_\sigma$} \rul{\forall \bm x\, \exists y^\sigma\, \vp(\bm x,y)}{\exists y^\sigma\, \forall\bm x\, \vp(\bm x,y)}{} \]
Same for $\izf$ and any other theory from Theorem \ref{sound}.
\end{theorem}
\begin{proof}
Use a reflexive pca over $\F$. WLOG, let $\bm x$ consist of a single variable. Suppose $\czf$ proves $\forall x\, \exists y^\sigma\, \vp(x,y)$. By soundness, let $\mb e$ be such that $\czf$ proves $\mb e\fo \forall x\, \exists y^\sigma\, \vp(x,y)$. According to our official convention, $\exists y^\sigma\, \vp(x,y)$ stands for 
\[ \forall z\, (\vartheta_\sigma(z)\imp \exists y\in z\, \vp(x,y)). \] 
Reasoning in $\czf$, we have that $\mb e\fo \exists y^\sigma\, \vp(\check x,y^\sigma)$ for every $x$. 
We know that there is a closed application term $\mb f$ such that $\mb f\fo \vartheta_\sigma(\dot F_\sigma)$. In particular, for every $x$ there exists $y\in F_\sigma$ such that $(\mb e\mb f)_0\simeq \bar y$ and $(\mb e\mb f)_1\fo \vp(\check x, \dot y)$. Note that $y$ does not depend on $x$. So let $y\in F_\sigma$ such that  $(\mb e\mb f)_0\simeq \bar y$.  We thus have that for every $x$, $(\mb e\mb f)_1\fo \vp(\check x, \dot y)$, and hence $\vp(x,y)$, as desired. 
\end{proof}

\subsection{Independence of premise}

\begin{theorem}[independence of premise rules]\label{premise rule}
Let $\psi(\bm x)$, $\vp(\bm x,y)$, $\theta(\bm x,\bm z)$ be formulas with displayed free variables.  Then $\czf$ is closed under the following rules:

\[ \tag{1} \rul{\forall \bm x\, (\neg \psi(\bm x)\imp \exists y^\sigma\, \vp(\bm x, y))}{ \exists y\, \forall \bm x\, (\neg \psi(\bm x)\imp y\in F_\sigma\land \vp(\bm x, y))}{} \]

\[ \tag{2} \rul{\forall \bm x\, (\neg \psi(\bm x)\imp \exists y^\sigma\, \vp(\bm x, y)) \qquad \qquad 
	\exists \bm x\, \neg\psi(\bm x)}{\exists y^\sigma\, \forall \bm x\, (\neg \psi(\bm x)\imp \vp(\bm x, y))}{} \]

\[ \tag{3} \rul{\forall \bm x\, (\forall \bm z\, \theta(\bm x,\bm z)\imp \exists y^\sigma\, \vp(\bm x, y)) \qquad \qquad  \forall \bm x\, \forall \bm z\, (\theta(\bm x, \bm z)\lor \neg \theta(\bm x, \bm z))}{\exists y^\sigma\, \forall \bm x\, (\forall \bm z\, \theta(\bm x,\bm z)\imp \vp(\bm x,y))}{} \]

\[ \tag{4} \rul{\forall \bm x\, (\forall \bm z^{\bm \rho}\, \theta(\bm x,\bm z)\imp \exists y^\sigma\, \vp(\bm x, y)) \qquad \qquad  \forall \bm x\, \forall \bm z^{\bm \rho}\, (\theta(\bm x, \bm z)\lor \neg \theta(\bm x, \bm z))}{\exists y\, \forall \bm x\, (\forall \bm z^{\bm \rho}\, \theta(\bm x,\bm z)\imp y\in F_\sigma\land \vp(\bm x,y))}{} \]
\smallskip 

Same for $\izf$ and any other theory from Theorem \ref{sound}.
\end{theorem}
\begin{proof}
For ease of notation, let $\bm x$ and $\bm z$ consist of a single variable $x$ and $z$ respectively.

Again, $\exists y^\sigma\, \vp(x,y)$ stands for  $\forall z\, (\vartheta_\sigma(z)\imp \exists y\in z\, \vp(x,y))$.
On the other hand, let $y\in F_\sigma\land \vp(x,y)$ be short for
\[  \forall z\, (\vartheta_\sigma(z)\imp y\in z\land \vp(x,y)). \]

(1) Use a reflexive pca over $\F$. Now suppose $\czf$ proves $\forall x\, (\neg \psi(x)\imp \exists y^\sigma\, \vp(x, y))$.
By soundness, let $\mb e$ be a closed application term such that $\czf$ proves $\mb e\fo \forall x\, (\neg \psi(x)\imp \exists y^\sigma\, \vp(x,y))$. From now one we argue in $\czf$. It follows from the definition of generic realizability that $\mb e\fo \neg \psi(\check x)\imp \exists y^\sigma\, \vp(\check x,y)$, for every $x$. 

For the sake of argument, suppose  that $\mb 0\fo \neg \psi(\check x)$. Then $\mb e\mb 0\fo \exists y^\sigma\, \vp(\check x,y)$. We know that there is a closed application term $\mb f$ such that $\mb f\fo \vartheta_\sigma(\dot F_\sigma)$. Then 
\[  \mb e\mb 0\mb f\fo \exists y\in \dot F_\sigma\, \vp(\check x,y).\]
Therefore $(\mb e\mb 0\mb f)_0\simeq \bar y$ for some (unique) $y\in F_\sigma$ and $(\mb e\mb 0\mb f)_1\fo \vp(\check x,\dot y)$. From this we conclude $\vp(x,y)$, since $\check x^\circ=x$ and $\dot y^\circ=y$.

It is now clear how to find $y$.  Simply, let 
\[ y= \{ u\in \bigcup F_\sigma \mid \exists v\in F_\sigma\, ( u\in v\land (\mb e\mb 0\mb f)_0\simeq \bar v)\}. \]

 We claim that $\forall x\, (\neg\psi(x)\imp y\in F_\sigma\land \vp(x,y))$. For this it is sufficient to note that $\neg \psi(x)$ implies  $\mb 0\fo \neg \psi(\check x)$. Note that, as $\dot F_\sigma$ is injective, the $y$ satisfying $\vp(x,y)$ is uniquely determined.

(2) Exercise for the reader. Use a reflexive pca over $\F$.

(3) Exercise for the reader. Apply $\sf UZR$ and ${\sf UR}_\sigma$. 

(4) follows from (1).

\end{proof}

\begin{remark}
The main use of reflexivity above consists in making sure that $\fo \vartheta_\sigma(\dot F_\sigma)$. In general, the argument goes through if for every type $\sigma$ there is a \emph{functional} $z\in \va$ such that $\fo \vartheta_\sigma(z)$, where by functional we mean that $\pair{a,y_1},\pair{a,y_2}\in z^*$ implies $y_1^\circ=y_2^\circ$. In this case, one can let
\[ y=\{u\in \bigcup F_\sigma\mid \exists a\in A\, \exists v\, ((\mb e\mb 0f)_0\simeq a\land u\in v^\circ\land \pair{a,v}\in z^*)\}, \]
where $f\fo \vartheta_\sigma(z)$. 
\end{remark}

\subsection{Explicit definability}

\begin{theorem}\label{def}
$\czf$ is closed under 

\[ \rul{\forall \bm x\, \exists y^\sigma\, \vp(\bm x, y)}{\exists ! y^\sigma\, (\delta(y)\land \forall \bm x\,  \vp(\bm x, y))}{\text{ for some formula } \delta(y)}  \]
\smallskip 

Same for $\izf$ and any other theory from Theorem \ref{sound}.
\end{theorem}
\begin{proof}
Use a  definable reflexive pca over $\F$. Note that in such case $\dot F_\sigma$ is also definable for every given type $\sigma$. An example of a pca over $\F$ definable in $\czf$ is given in the upcoming and final section. 
\end{proof}

\section{A direct construction in $\czf$ of a reflexive pca over finite types} \label{example}

We first describe the general idea in a classical  setting. Recall that $\F=\bigcup \Ff$ with $\Ff=\{F_\sigma\mid \sigma\in\ft\}$.  Let function application  be given by  $f x\simeq y$ iff $f$ is a function, $x\in\dom(f)$ and $f(x)=y$. Since $\F$ is closed under function application, this gives us a partial algebra on $\F$.  The  idea would be to define a partial application map  on the powerset $\ps(\F)$ and take $x\mapsto \{x\}$ to be the embedding. Unfortunately, the usual constructions on $\ps(\F)$ do not yield a pca.\footnote{We remind the reader of a construction due to van Oosten. Given any pca $A$, one can define a partial binary operation on  $\ps(A)$   by letting $XY\simeq Z$ iff $Z=\{ab\mid a\in X \land b\in Y\}$ and application is total on $X\times Y$, that is, $ab$ is defined for all $a\in X$ and $b\in Y$. This need not be a pca with combinators $\{\mb k\}$ and $\{\mb s\}$. Note that in general $\{\mb s\}XYZ$ is smaller than $XZ(YZ)$. However, the totality requirement makes it an ordered pca in the sense of van Oosten \cite[1.8]{Oosten08}.  A similar contruction on $\F$ gives rise to an ordered pca on $\ps(\F)$, but this fails to be a pca for the same reasons.} To solve this, we  introduce the notion of arity and  work with nonempty subsets of $\F\times \ar$, where $\ar$ is the set of arities.   We use arities, that is types of the form $p_0\cdots p_n\to q$, to iterate function application in a prescribed  manner. For example, we can assign  arity $\sigma\tau\to\rho$ to a function  $f$ of type $\sigma\to\tau\to\rho$ and thus see $(f,\sigma\tau\to\rho)$ as a function from $F_\sigma\times F_\tau$ to $F_\rho$ by currying. We then define a partial application map on nonempty subsets of $\F\times \ar$, so that the resulting partial algebra is a pca that embeds $\F$ via the canonical embedding $x^\sigma\mapsto\{(x,\sigma)\}$. We require the sets to be nonempty so that the combinator $\mb k$ does its job. However, in order to obtain a reflexive pca, in particular a  \lq definition by cases\rq\ combinator,   we also have to enlarge  the type structure  $\Ff$ by allowing (enough) dependent products. Formally, this is how we proceed. 

Let $\Gf$ be inductively defined by:
\begin{itemize}
\item $\omega\in \Gf$
\item if $F,G\in \Gf$, then $F\to G\in \Gf$;
\item if $F_n\in \Gf$ for every $n\in\omega$, then $\prod_{n\in \omega}F_n\in \Gf$,
\end{itemize}
where in general
\[  \prod_{x\in F}G_x=\{ f\colon F\to \bigcup_{x\in F}G_x\mid \forall x\in F\, (f(x)\in G_x)\}. \]

Set $\Gbb=\bigcup \Gf$. Note that $\F\subseteq \Gbb$. On elements of $\Gbb$ we will always consider function application.
We define the set of arities $\ar$ by the following inductive clauses:
\begin{itemize}
	\item $o\in\ar$;
	\item if $p_0,\ldots,p_n,q\in\ar$ then $p_0\cdots p_n\to q\in \ar$;
	\item if $p_n\in\ar$ for every $n\in\omega$, then $\prod_{n\in\omega}p_n\in \ar$.
\end{itemize}
Note that $\ft\subseteq\ar$. If $\sigma,\tau\in\ft$, we denote by $\sigma^n\to\tau\in P$ the arity
\[ \overbrace{\sigma\cdots\sigma}^{n+1}\to \tau. \]

Let $\Gr=\Gbb\times \ar$. We use $\mt x,\mt y,\mt z,\ldots$ to denote elements of $\Gr$. We define a partial function from $\Gr\times \Gr^{<\omega}$ to $\Gr$ by letting $\mt x(\vec{\mt y})\simeq \mt z$
iff $\mt x=(x,p)$, $\mt z=(z,q)$ and  either one of the following applies:
\begin{itemize}
	\item[(i)] $p=p_0\cdots p_n\to q$, $\vec{\mt y}=\langle (y_i,p_i)\mid i\leq n\rangle$,  and $xy_0\cdots y_n\simeq z$;
	\item[(ii)] $n=0$, $p=\prod_{n\in\omega}p_n$,  $\vec{\mt y}=\langle (m,o)\rangle$ for some $m\in\omega$, $p_m=q$ and $xm\simeq z$.
\end{itemize}
The reader should keep in mind that $xy_0\cdots y_n$ and $xm$  are defined by (iterated) function application. We usually write $\mt x(\mt y_0,\ldots,\mt y_n)$ for $\mt x(\langle \mt y_0,\ldots,\mt y_n\rangle)$.

Let $\Gr^*=\ps(\Gr)\setminus \{\emptyset\}$. A partial application map on $\Gr^*$ is then defined by letting $ab\simeq c$ iff
\[ c=\{ \mt z\in \Gr\mid \exists \mt x\in a\, \exists \mt y_0,\ldots,\mt y_n\in b\, (\mt x(\mt y_0,\ldots,\mt y_n)\simeq \mt z)\}. \]
We now equip  $\Gr^*$ with combinators. \medskip

(1) Let us define $\mb k\in\Gr^*$.  First, for all $F,G\in \Gf$, let $k_{FG}\colon F\to G\to F$ be the unique function such that $kx^Fy^G=x$.  Note that all $k$'s  are in $\Gbb$ since $\Gf$ is closed under exponentiation.
Let 
\[  \mb k=\{ (k_{FG},p\to q\to p)\mid  F,G\in \Gf\land p,q\in\ar\}. \]
{}\smallskip

(2) We define $\mb s\in\Gr^*$ as follows.  
For any choice of $n,m,n_0,\ldots,n_m\in\omega$ we consider all  functions $s$'s in $\Gbb$  such that 
\begin{multline*}   sxy_0\cdots y_mz_0\cdots z_nz_{00}\cdots z_{0n_0}\cdots\cdots z_{n0} \cdots z_{mn_m}= \\ =xz_0\cdots z_n(y_0z_{00}\cdots z_{0n_0})\cdots (y_mz_{m0} \cdots z_{mn_m}). 
\end{multline*}
We suppress the type information for notational convenience. Note  that all  $s$'s are in $\Gbb$ since $\Gf$ is closed under exponentiation.  Any such function can be assigned an arity of the form 
\[ p\to r_0\cdots r_m\to q_0\cdots q_nq_{00}\cdots q_{0n_0}\cdots\cdots q_{m0} \cdots q_{mn_m}\to q, \]
where $p=q_0\cdots q_n\to p_0\cdots p_m\to q$ and  $r_i=q_{i0}\cdots q_{in_i}\to p_i$ for every $i\leq m$.

Let
\[ \mb s=S\cup C, \]
where $S$ consists of all pairs $(s,r)$, where $s$ is as above and the arity $r$ agrees with the type of $s$ in  the sense just described, and 
\[ C= \{ (\lambda x^F y^G.0^o,p\to q\to o)\mid  F,G\in \Gf\land p,q\in\ar\}. \]
We add $C$ only to ensure that $\mb sab\downarrow$ for all $a$ and $b$.\medskip 

(3) Numerical combinators are easily definable. Let $\succe=\{(\lambda n.n+1,o\to o)\}$ and $\pred=\{(\lambda n.n-1,o\to o)\}$.\medskip

(4) As for $\mb d\in\Gr^*$, for all $F,G\in \Gf$, let $d_{FG}\colon \prod_{n\in\omega}\prod_{m\in\omega}F_{nm}$, where 
\[ F_{nm}=\begin{cases}
	F\to G\to F & \text{ if $n=m$;} \\
	F\to G\to G & \text{ otherwise,}
\end{cases} \]
such that 
\[ d_{FG}nm\simeq\begin{cases}
	k_{FG} & \text{ if } n=m; \\
	\bar k_{FG} &  \text{ otherwise. }
\end{cases}\]
Here, $\bar k_{FG}\colon F\to G\to G$ is defined by $\bar k x^Fy^G=y$.  It is easy to check that for every $F,G\in \Gf$, 
\[ \prod_{n\in\omega}\prod_{m\in\omega}F_{nm}\in \Gf, \]
and therefore $d_{FG}\in \Gbb$. Similarly, for all $p,q\in\ar$, let $(p,q)\in \ar$ be defined as $\prod_{n\in\omega}\prod_{m\in\omega}p_{nm}$, where
\[ p_{nm}=\begin{cases}  p\to q\to p & \text{ if $n=m$;} \\
	p\to q\to q  & \text{ otherwise.}
	\end{cases} \] 
Let
\[  \mb d=\{ (d_{FG},(p,q))\mid  F,G\in \Gf\land p,q\in\ar \}. \]
{}\smallskip

(5) Combinator   $\mb i_{\sigma\tau}\in\Gr^*$ for  $\sigma,\tau\in\ft$. Let
\[ \mb i_{\sigma\tau}=\{(i_{\sigma\tau}^n,(\sigma^n\to \tau)\to \sigma\to\tau)\mid n\in\omega\},\]
where $i_{\sigma\tau}^n\colon (\overbrace{F_\sigma\to\cdots \to F_\sigma}^{n+1}\to F_\tau)\to F_\sigma\to F_\tau$ is defined by 
$i^n_{\sigma\tau}fx=f \overbrace{x\cdots x}^{n+1}$.\medskip

One could show  in $\izf$ that $\Gr^*$ is a reflexive pca over $\F$. On the other hand, $\Gr^*$ is not even a set in $\czf$. Note that in $\czf$ we could inductively define $\Gbb$ as a class. For our purposes however, a class will not do. It turns out that we just need very few dependent products, as can be easily gleaned from the construction above.  Also, we have to settle on the right notion of nonemptiness. But this is easily arranged: we just take inhabited sets ($x$ is inhabited if $\exists u\, (u\in x)$). The construction in $\czf$ proceeds as follows.

First, we can form the set   $\Gf=\bigcup_{s\in\omega}\Gf_s$,  where
\begin{itemize}
	\item $\Gf_0=\{\omega\}$;
	\item $\Gf_{s+1}=\Gf_s\cup\{ F\to G\mid F,G\in \Gf_s\}\cup \{ \prod_{n\in\omega}F_s\mid \forall n\in\omega\, (F_n\in \Gf_s)\}$. 
\end{itemize}
Now, $\omega\in \Gf$ and $\Gf$ is closed under exponentiation. Let  $\Gbb=\bigcup \Gf$. It follows that $F_\sigma\in \Gf$ for every finite type $\sigma$, and so  $\F\subseteq \Gbb$. In a similar manner, we obtain in $\czf$ a sufficiently large set $P$ of arities.   Let $P=\bigcup_{s\in\omega}P_s$,  where
\begin{itemize}
	\item $P_0=\{o\}$;
	\item $P_{s+1}=P_s\cup\{ p_0\cdots p_n\to q \mid p_0,\ldots,p_n,q\in P_s\}\cup \{ \prod_{n\in\omega}p_n\mid \forall n\in\omega\, (p_n\in P_s)\}$. 
\end{itemize}

We can then form the set $\Gr=\Gbb\times \ar$.  Finally, unless we work in $\izf$, where we have access to the full power set, we need to get by with sufficiently (set-)many  subsets of $\Gr$. This is how we can proceed. Consider 
\[ X=\{\mb k,\mb s, \mb d\}\cup\{\{(x,\sigma)\}\mid  \sigma\in\ft, x\in F_\sigma\}\cup\{\mb i_{\sigma\tau}\mid  \sigma,\tau\in\ft\}. \]
We can construct the total algebra $\bar X\subseteq \Gr^*$  generated by $X$ under application in  $\Gr^*$. To see this, we define by recursion
\begin{align*}
X_0 & = X \\
X_{s+1} & = X_s \cup \{ab\mid a,b \in X_s\} 
\end{align*}
and set $\bar X=\bigcup_{s\in\omega} X_s$. 
 We then let 
\[ A=\{a\in \bar X\mid  a \text{ is inhabited}\}.  \] 
Note that $A$ consists of inhabited subsets of $\Gr$. Also, for all $a,b\in A$, we have $ab\downarrow$ iff $ab$ is inhabited iff $ab\in A$.

\begin{theorem}[$\czf$]
$A$ is a  reflexive pca over $\F$ with embedding $x^\sigma\mapsto \bar x=\{(x,\sigma)\}$.
\end{theorem} 
\begin{proof}
First, the function $x\mapsto \bar x$ from $\F$ to $A$ is given by 
\[ \{\langle x, (x,\sigma)\rangle\mid x\in F_\sigma\land \sigma\in\ft\}. \]
By induction on the type, one can verify that if $F_\sigma\cap F_\tau$ is inhabited, then $\sigma=\tau$. 
Therefore the set above is indeed a function.  That this map provides an embedding of partial algebras from $\F$ into $A$  is  immediate.  In fact, if $f\in F_{\sigma\tau}$ and $x\in F_\sigma$, then  $(f,\sigma\tau)(x,\sigma)\simeq (f(x),\tau)$.

Let us check the combinators. 

(1) Combinator $\mb k$. Let $a,b\in A$. We have
	\[   \mb kab=\bigcup_{F,G\in \Gf}\{(k_{FG}x^Fy^G,p)\mid  (x,p)\in a \land \exists q\, ((y,q)\in b)\}=a. \]
Use the fact that  $a$ and $b$ are inhabited. 
	
(2) Combinator $\mb s$. Notice that $\mb sab\downarrow$ since $(0,o)\in \mb s a b$. Now, $\mb sabc$ is designed to contain exactly  all elements of $ab(cb)$. The verification that $\mb sabc\simeq ac(bc)$ is  a simple exercise. 
	
(3)  The verification that $\succe$ and $\pred$ behave as desired is immediate. 
	
(4) Combinator $\mb d$. It is not difficult to see that 
	\[ \mb d\bar n\bar m=\begin{cases}
		\mb k & \text{ if  $n=m$}; \\
		\overline{\mb k} & \text{ otherwise;}
	\end{cases} \] 
	where $\overline{\mb k}=\{(\bar k_{FG},p\to q\to q)\mid  F,G\in \Gf \land p,q\in\ar\}$, and $\bar k_{FG}=\lambda x^Fy^G.y$. Therefore $\mb d\bar n\bar m ab\simeq a$ if $n=m$ and $b$ otherwise.
	
(5)  Combinator $\mb i_{\sigma\tau}$. Let $F=F_\sigma$ and $G=F_\tau$. Suppose that $a\in A$ and that for every $x\in F$ there exists $y\in G$ such that $a\bar x\simeq \bar y$. We want to show that $\mb i_{\sigma\tau}a\simeq \bar f$, where $f\colon F\to G$ is such that 
\[ \tag{$*$} a\bar x\simeq \overline{f(x)}  \ \text{ for every $x\in F$}.\]
Recall that $\mb i_{\sigma\tau}=\{(i_{\sigma\tau}^n,(\sigma^n\to \tau)\to \sigma\to\tau)\mid n\in\omega\}$.
Let us denote by $\mt i^n_{\sigma\tau}$ the  $n$-th element of  $\mb i_{\sigma\tau}$.  Let us also write $F\to^n G$ for 
\[   \overbrace{F\to\cdots \to F}^{n+1}\to G.\]
Note the difference between the set $F\to^n G$ and the arity $\sigma^n\to\tau$.
Recall that \[i_{\sigma\tau}^n \colon (F\to^n G)\to F\to G. \] 
By definition, 
$\mb i_{\sigma\tau}a=\{\mt i^n_{\sigma\tau}(\mt g)\mid n\in\omega\land \mt g\in a\}$.

Let us point out that $\Gf$ is sparse in the sense that if $F,G\in \Gf$ have some overlap, meaning that $F\cap G$ is inhabited, then $F=G$. This is proved by induction on $F\in \Gf$ (i.e., by induction on the stage $n\in\omega$ such that $F\in \Gf_n$). This feature is crucial here and will be used without further notice. 

We start off by proving that $\bar f\subseteq \mb i_{\sigma\tau}a$, that is,  $(f,\sigma\tau)\in \mb i_{\sigma\tau}a$.  Pick any $u\in F$. We can do this since every $F\in\Ff$ is inhabited.  By $(*)$,  there must be a $\mt g=(g,p) \in a$ such that 
\[  \mt g(\overbrace{(u,\sigma),\cdots,( u,\sigma)}^{n+1})\simeq (f(u),\tau). \]
It is not too difficult to check that $g\colon  F^n\to G$ and $\mt g(\langle (x,\sigma)\mid i \leq n\rangle)\downarrow$  for every $x\in F$.
Notice that $p=\sigma^n\to\tau$ or $p=\prod_{m\in\omega}p_m$, in which case $n=0$, $\sigma=o$ and $p_m=\tau$ for every $m\in\omega$.  
Since $\mt g(\langle (x,\sigma)\mid i\leq n\rangle) \in a\bar x$, it thus follows by $(*)$ that  $\mt g(\langle (x,\sigma)\mid i\leq n\rangle)\simeq (f(x),\tau)$ for every $x\in F$. But then $\mt i^n_{\sigma\tau}(\mt g)\simeq(f,\sigma\tau)$. This shows one direction.    

The other direction $\mb i_{\sigma\tau}a\subseteq \bar f$ is similar. Let  $\mt g=(g,p)\in a$ and suppose that $\mt i^n_{\sigma\tau}(\mt g)\downarrow$, so that $\mt i^n_{\sigma\tau}(\mt g)\in \mb i_{\sigma\tau}a$.  Any element of $\mb i_{\sigma\tau}a$ is obtained this way.  By definition, it must be $p=\sigma^n\to \tau$ and $i_{\sigma\tau}^ng\downarrow$. In particular,  $g\colon F\to^n G$.  It thus follows that for every $x\in F$ there exists  $y_x\in G$ such that $\mt g(\langle (x,\sigma)\mid i\leq n\rangle)\simeq (y_x,\tau)$.
On the other hand, for every $x\in F$, $(y_x,\tau)\in a\bar x=\overline{f(x)}$, and so $y_x=f(x)$. Therefore  $\mt i^n_{\sigma\tau}(\mt g)\simeq (f,\sigma\tau)$, as desired. 
\end{proof}

\end{document}